\documentclass[10pt]{amsart}
\usepackage{amssymb}

\newtheorem{theorem}{Theorem}
\newtheorem{hypothesis}[theorem]{Hypothesis}
\newtheorem{lemma}[theorem]{Lemma}

\newtheorem{corollary}[theorem]{Corollary}

\newcommand{\R}{{\mathbb R}}

\newcommand{\dstar}{d^*\!}

\newcommand{\xnot}{x_0}
\renewcommand{\l}{\lambda}

\newcommand{\demo}{\noindent{\it Proof.} }
\newcommand{\findemo}{\qed\bigskip}

\begin{document}

\title[Spectral cluster bounds]
{Subcritical $L^p$ bounds on spectral clusters\\
for Lipschitz metrics}
\thanks{The authors were supported in part by NSF grants DMS-0140499,
DMS-0354668, DMS-0301122, and DMS-0354539.}
\author{Herbert Koch}
\address{Mathematisches Institut, Universit\"at Bonn, 53113 Bonn, Germany}
\email{koch@math.uni-bonn.de}
\author{Hart F. Smith}
\address{Department of Mathematics, University of Washington, 
Seattle, WA 98195}
\email{hart@math.washington.edu}
\author{Daniel Tataru}
\address{Department of Mathematics, University of California,
Berkeley, CA 94720}
\email{tataru@math.berkeley.edu}

\begin{abstract}
We establish asymptotic bounds on the $L^p$ norms of
spectrally localized functions in the case of 
two-dimensional Dirichlet forms with
coefficients of Lipschitz regularity. 
These bounds are new for the range $6<p<\infty$.
A key step in the proof is bounding the rate at which energy spreads
for solutions to hyperbolic equations with Lipschitz coefficients.

\end{abstract}

\maketitle

\section{Introduction}
The purpose of this paper is to establish $L^p$ bounds on eigenfunctions,
or more generally spectrally localized functions, associated to 
Dirichlet forms on a compact manifold. 
The question of interest is the dependance of the bounds on the
H\"older regularity of the coefficients of the form.
We consider here the case of
Dirichlet forms with Lipschitz coefficients for simplicity, but the proofs
can be adapted to the case of $C^s$ coefficients, where $0<s<2$.
Our work is
restricted to the case of two-dimensional manifolds, however.

Consider the eigenvalue problem for
a Dirichlet form, where we work on a compact manifold $M$ without boundary,
\begin{equation*}%\label{quadform}
\dstar(a\,d\phi)+\l^2\rho\,\phi=0\,.
\end{equation*}
Here, $a$ is a section of real,
symmetric quadratic forms on $T^*(M)$, with associated linear transforms
$a_x:T^*_x(M)\rightarrow T_x(M)$,
and $\rho$ is a real valued function on $M$. Here, $\dstar$ denotes the adjoint
of $d$ relative to a fixed volume form $dx$.
We assume both $a$ and $\rho$ are strictly positive, with
uniform bounds above and below.
We note that this setting includes the Laplace-Beltrami operator
on a Riemannian manifold. The parameter $\l\ge 0$ is referred to as the
frequency of the eigenfunction $\phi$.

A {\it spectral cluster} of frequency $\l$ 
is a combination of eigenfunctions with frequencies in the range
$[\l-1,\l]$. In the case that $a$ and $\rho$ are smooth,
Sogge \cite{So} established the following best possible 
$L^p$ bounds on spectral clusters,
\begin{equation}\label{Soggebound}
\bigl\|f\bigr\|_{L^p(M)}\lesssim
\begin{cases}\,\l^{\frac{n-1}2(\frac 12-\frac 1p)}\,\|f\|_{L^2(M)}\,,
\quad &2\le p\le p_n\\
\l^{n(\frac 12-\frac 1p)-\frac 12}\,\|f\|_{L^2(M)}\,,
&p_n\le p\le\infty
\end{cases}
\end{equation}
The critical index is $p_n=\frac {2(n+1)}{n-1}$.
Semiclassical generalizations were obtained by 
Koch-Tataru-Zworksi \cite{KTZ}.
The bounds \eqref{Soggebound}
hold in case $a$ and $\rho$ are of regularity $C^{1,1}$
by \cite{Sm1},
but based on an observation of Grieser \cite{Gr}, and examples of Smith-Sogge
\cite{SmSo1} and the authors \cite{KST}, 
they fail for coefficients of $C^s$ regularity if $s<2$.

For metrics of regularity $C^s$ with $s<2$ (or Lipschitz in case $s=1$)
best possible $L^p$ bounds on spectral clusters have been established
on the range $2\le p\le p_n$, as well as for $p=\infty$; see \cite{Sm2}
for the case $1\le s<2$, and \cite{KST} for the case $s<1$. 
This leaves open the {\it subcritical} case $p_n<p<\infty$, where
the upper bounds on the exponent of $\l$
that can be obtained from \cite{KST} and \cite{Sm2}
by interpolation do not match the lower bounds
that follow from the examples of \cite{KST} and \cite{SmSo1}.

In this paper we obtain bounds for
$p_n<p<\infty$, for Lipschitz coefficients and $n=2$, 
which improve upon the results of \cite{Sm2}. 
They do not match the exponent displayed by the Rayleigh whispering
mode example noted in \cite{Gr}, but the
difference is exponentially small as $p\rightarrow\infty$.
Our results are restricted to $n=2$,
but all steps adapt to $C^s$ coefficients for $0<s<2$, and improve
upon \cite{KST} and \cite{Sm2} on this range of $p$.

Thus, consider a Dirichlet form on a two-dimensional compact manifold without
boundary, with $a$ and $\rho$ of Lipschitz regularity. Let
$$
\gamma(p)=2\Bigl(\tfrac 12-\tfrac 1p\Bigr)-\tfrac 12
$$
be the exponent occuring in the subcritical estimates \eqref{Soggebound}.

By  Theorem 2 of \cite{Sm2}, in this case {\it no-loss estimates} hold
on cubes $Q$ of sidelength $\l^{-\frac 13}$,
\begin{equation}\label{noloss}
\|f\|_{L^p(Q)}\lesssim\l^{\gamma(p)}\|f\|_{L^2(M)}\,,\qquad 6\le p\le\infty\,.
\end{equation}
The Rayleigh whispering mode examples show that
if $p=6$ the size of $Q$ cannot be
increased without increasing the exponent. The main result of this
paper is that, for larger $p$, the following {\it log-loss estimates} hold on
cubes $Q$ of sidelength $\l^{-\frac 13 2^{(6-p)/2}}$
\begin{equation}\label{logloss}
\|f\|_{L^p(Q)}\lesssim (\log\l)^{p-6}\l^{\gamma(p)}\|f\|_{L^2(M)}\,,
\qquad p=6,8,10,12,\ldots
\end{equation}
If $f$ is conically microlocalized in frequency, then 
(see \cite[(14)-(15)]{Sm2})
$Q$ can be replaced
by a thin slab of size $1\times \l^{-\frac 13 2^{(6-p)/2}}$. Summing
over such slabs, one obtains the following.

\begin{theorem}\label{theorem1}
Suppose that $a$ and $\rho$ are of Lipschitz regularity, on a two-dimensional
compact manifold without boundary. Then, for $p=6,8,10,12,\ldots$
\begin{equation}\label{KSTbound}
\|f\|_{L^p(M)}\lesssim (\log\l)^{p-6}\l^{\sigma(p)}\|f\|_{L^2(M)}\,,
\qquad
\sigma(p)=\gamma(p)+\tfrac 1{3p}\,2^{\frac{6-p}2}\,.
\end{equation}
\end{theorem}

\medskip

To place this result in context, the Rayleigh mode and a reflection
argument shows that, for Lipschitz coefficients and $n=2$,
one cannot establish better estimates
than the following
\begin{equation}\label{examples}
\bigl\|f\bigr\|_{L^p(M)}\lesssim
\begin{cases}\,\l^{\frac 23(\frac 12-\frac 1p)}\,\|f\|_{L^2(M)}\,,
\quad &2\le p\le 8\\
\l^{\gamma(p)}\,\|f\|_{L^2(M)}\,,
&8\le p\le\infty
\end{cases}
\end{equation}
The results of \cite{Sm2} imply that the bounds \eqref{examples} hold
for general Lipschitz coefficients for the cases $2\le p\le 6$ and $p=\infty$.
The exponent $\sigma(p)$
in \eqref{KSTbound}, which agrees with that of \eqref{examples} 
for $p=6$, misses \eqref{examples} for $6<p<\infty$ by a
factor which decays exponentially as $p\rightarrow\infty$.

We remark that
the bounds \eqref{examples} were established in \cite{SmSo2}
for smooth Dirichlet forms on two-dimensional manifolds with
boundary, with either Dirichlet or Neumann conditions at the boundary.
A manifold with boundary can be thought of as a special case of
a Lipschitz metric, by reflecting coefficients normally across the boundary. 
The example of \cite{SmSo1} for Lipschitz metrics and $n=2$
is generated by reflecting a Rayleigh whispering mode from the unit 
disc.

The proof of \eqref{logloss} is inductive. The estimate for $p+2$
is derived from the estimate for $p$, together with an almost
orthogonal decomposition of $f$ into tubular pieces. 
Essentially, one can localize $f$ in frequency to a cone of angle $\delta$, 
and in space to a tube of diameter $\delta^2$,
and control the energy flow over distance $\delta$. 
For this reason, the diameter of the log-loss cubes for $p+2$
is the square root of the diameter of the log-loss cubes for $p$. 
The argument that allows summation over different tubes with
a $(\log\l)^2$ loss works only for $n=2$, however. 
Improving Theorem \ref{theorem1}
appears then to hinge on controlling energy flow over longer distances,
and improving the summation argument to allow $n\ge 3$.

The bounds we establish hold more generally 
for functions that satisfy a quasimode condition
\begin{equation}\label{quasimode}
\dstar(a\,df)+\l^2\rho\,f=\dstar g_1+g_2\,.
\end{equation}
If $f$ is a spectral cluster, then \eqref{quasimode} holds on $M$
with $g_1=0$ and $\|g_2\|_{L^2}\lesssim \l\|f\|_{L^2}$.
Allowing the term $g_1$ makes localization arguments simpler.
In particular \eqref{quasimode} holds, with
$\|g_1\|_{L^2}\lesssim \|f\|_{L^2}$ and 
$\|g_2\|_{L^2}\lesssim \l\|f\|_{L^2}$,
for the product of a spectral cluster $f$ with a unit
size bump function, so we may assume that we work in local coordinates.
After rescaling and extending, we may assume that $a$ and $\rho$ are 
globally close on $\R^2$ to the flat metric, 
$$%\begin{equation}\label{lipcond}
\|a^{ij}-\delta^{ij}\|_{Lip(\R^2)}+
\|\rho-1\|_{Lip(\R^2)}\le c_0\,.
$$%\end{equation}
We establish the estimate \eqref{logloss} by an induction argument,
for which the starting point is the localized version of
the no-loss estimate \eqref{noloss}
for $p=6$. At each step of the induction
$p$ increases by 2, and we establish estimates on
cubes of square-root the sidelength of the previous step. 
A loss of $(\log\l)^2$ is incurred at each step, however.
The hypothesis and induction argument are as follows.

\begin{hypothesis}\label{hypo}
Suppose that equation \eqref{quasimode} holds on $Q^*$,
where $Q^*$ denotes the double of the cube $Q$.
Then the following inequality holds, where
$\ell(Q)$ denotes the sidelength of $Q$
\begin{multline*}
\|f\|_{L^p(Q)}\le C_p\,
(\log\l)^{p-6}\l^{\gamma(p)}\Bigl(\,
\ell(Q)^{-\frac 12}\|f\|_{L^2(Q^*)}
+\l^{-1}\ell(Q)^{-\frac 12}\|df\|_{L^2(Q^*)}\\
+\ell(Q)^{\frac 12}\,\|g_1\|_{L^2(Q^*)}
+\l^{-1}\ell(Q)^{\frac 12}\|g_2\|_{L^2(Q^*)}
\Bigr)\,.
\end{multline*}
\end{hypothesis}

\begin{theorem}
Assume
that Hypothesis \ref{hypo} holds for a given $p\in [6,\infty)$,
uniformly over cubes $Q$ of a given sidelength $\ell(Q)=\delta^2$, where
$1\ge \delta \ge \l^{-\frac 16}$. 
Then Hypothesis \ref{hypo} holds with $p$ replaced by $p+2$,
uniformly over cubes $Q$ of sidelength $\ell(Q)=\delta$.
\end{theorem}

We remark that the norm on the right hand side in Hypothesis \ref{hypo}
should be thought of
as the $L^2$-energy of $f$ on $Q^*$. If the functions involved are
localized to frequencies $\xi$ of magnitude $\l$, and $\xi$ in a small
cone about the $\xi_1$ axis, then the right hand side is a
replacement for
$\|f\|_{L^\infty_{x_1}L^2_{x_2}}+\|Pf\|_{L^1_{x_1}L^\infty_{x_2}}$.

The outline of this paper is as follows. In Section \ref{section:tube},
we establish the key decomposition of $f$ as a sum of terms, each supported
in a thin geodesic tube of dimensions $\delta\times\delta^2$.
This decomposition at multiple scales is inspired by the work of Geba-Tataru
\cite{GeTa}.
In Section \ref{section:overlap} we establish
$\ell^q$ bounds on the overlaps of collections of such tubes,
which is applied in Section \ref{section:proof} to complete the proof.

For the remainder of this section we carry out some simple
initial reductions.
Consider a cube $Q_0$ of sidelength $\delta$. Let $\psi$ be a scaled
bump function, supported in $Q_0^*$ and equal to $1$ on $Q_0$.
Then Hypothesis \ref{hypo}
with $\ell(Q)=\delta$ is unchanged if we replace $f$ by
$\psi f$, hence we may assume $f$ is supported in $Q_0^*$,
and in particular use $\|f\|_{L^2}$ instead of $\|f\|_{L^2(Q_0^*)}$.

Next, we split $f$ into components $f=f_{<\l}+f_\l+f_{>\l}$, by localizing
respectively to
frequencies smaller than $c^2\l$, comparable to $\l$, and larger
than $c^{-2}\l$, where $c$ is a fixed
small constant.
By the arguments of \cite[Corollary 5]{Sm2},
$$
\l\,\|f_{<\l}\|_{L^2}+\|df_{>\l}\|_{L^2}\lesssim
\|f\|_{L^2}
+\l^{-1}\|df\|_{L^2}
+\|g_1\|_{L^2}
+\l^{-1}\|g_2\|_{L^2}\,.
$$
Since $\l^{2(\frac 12-\frac 1p)-1}\le\l^{\gamma(p)}\ell(Q_0)^{\frac 12}$,
Sobolev embedding yields that Hypothesis \ref{hypo}
holds with $\|f_{<\l}\|_{L^p}$ and $\|f_{>\l}\|_{L^p}$ on the left hand side.
Thus we restrict attention to the case that $f$ is frequency localized
to $|\xi|\approx\l$. By writing $f_\l$ as a finite sum of terms,
we may also assume that $f$ is frequency localized to $|\xi_2|\le c\l$.

Define the operator 
$$
P_{\l\delta}=\dstar a_{\l\delta}\,d+\l^2\rho_{\l\delta}\,,
$$
where the coefficients $a_{\l\delta}$ and $\rho_{\l\delta}$
are smoothly truncated in frequency to $|\xi|\le c\l\delta$.

Provided $\ell(Q)\le \delta$, then
Hypothesis \ref{hypo} is unchanged if we replace the defining equation by
$P_{\l\delta}f=\dstar g_1+g_2$, since the difference $(P-P_{\l\delta})f$
can be absorbed into $g_1$ and $g_2$, leaving the right hand side
of the inequality unchanged up to a constant.

Given a cube $Q$ and parameters $\delta$, $\l$, we set
$$
|||f|||_{\l,\delta,Q}=
\delta^{-\frac 12}\|f\|_{L^2(Q)}
+\l^{-1}\delta^{-\frac 12}\|df\|_{L^2(Q)}
+\l^{-1}\delta^{\frac 12}\|P_{\l\delta}f\|_{L^2(Q)}\,.
$$
We use $|||f|||_{\l,\delta}$ to denote the norm in case $Q=\R^2$.

Since $f$ is frequency localized to $|\xi|\approx\l$, as is
$P_{\l\delta}f$, we may absorb the term $\dstar g_1$ into $g_2$.
Thus, by the preceeding comments, we are reduced to the following.

\begin{theorem}\label{maintheorem'}
Suppose that $f$ is frequency localized to $|\xi|\approx \l$ and 
$|\xi_2|\le c\l$.
Then the following holds, uniformly on cubes $Q$ of sidelength $\delta$,
$$
\|f\|_{L^{p+2}(Q)}\lesssim
(\log \l)^{p-4}\l^{\gamma(p+2)}|||f|||_{\l,\delta}
$$
under the assumption that the following holds, uniformly on cubes
$Q$ of sidelength $\delta^2$,
$$
\|f\|_{L^p(Q)}\lesssim
(\log \l)^{p-6}\l^{\gamma(p)}|||f|||_{\l,\delta^2,Q^*}
$$
\end{theorem}

\section{The tube decomposition}\label{section:tube}
Let $f$ be as in Theorem \ref{maintheorem'}, and fix a cube $Q_0$
of sidelength $\delta$ and center $\xnot$. As above, let $\psi=1$ on
$Q_0$ and vanish outside $Q_0^*$. In this section we produce a decomposition
\begin{equation}\label{decomp}
\psi f=\sum_{T \in \mathcal T}f_T+f_0\,,
\end{equation}
where $f_0$ is an error term whose $L^p$ norms can be 
appropriately bounded by Sobolev embedding.
Each $f_T$ is compactly supported in a tube $T$.
The index $T$ varies over a collection $\mathcal T$ of
tubes of diameter $\delta^2$ and length $\delta$,
associated to bicharacteristic directions
of $P_{\l\delta}$ at angular separation $\delta$. Each $\widehat f_T$ is
concentrated (in a weighted $L^2$ sense)
in a ball of diameter $\l\delta$. We further have the bounds,
\begin{equation}\label{sqsum}
\Biggl(\,\sum_T 
|||f_T|||^2_{\l,\delta}\,\Biggr)^{\frac 12}
\le C \, |||f|||_{\l,\delta}\,.
\end{equation}

Let $\Gamma$ be the characteristic set of $P_{\l\delta}$ which
lies near $Q_0^*\times\text{support}(\hat f)$,
$$
\Gamma=\bigl\{(x,\xi):
\langle a_{\l\delta}(x)\,\xi,\xi\rangle=\l^2\rho_{\l\delta}(x)\bigr\}
\cap 
\,Q_0^*\times\{|\xi_2|\le 2c\l\}\,.
$$
Since $a_{\l\delta}$ and $\rho_{\l\delta}$ 
are pointwise close to the flat metric,
the set $\Gamma$ can be realized as the union of two graphs 
$\xi_1=\gamma_\pm(\xi_2)$. Since
$a_{\l\delta}$ and $\rho_{\l\delta}$ are Lipschitz 
and $|x-\xnot|\le \delta$, the characteristic
set $\Gamma$ is contained in a $\l\delta$ neighborhood of $\Gamma_{\xnot}$.
Let $q_j(\xi)$ be a finite-overlap cover of $\Gamma$ by
$\approx \delta^{-1}$ smooth bump functions, each supported in a 
ball of diameter $\approx\l\delta$ centered on $\Gamma_{\xnot}$,
so that $\phi(\xi)=1-\sum_j q_j(\xi)$ vanishes on a $\l\delta$ size 
neighborhood of $\Gamma$. Thus, 
$\phi(\xi)P_{\l\delta}(x,\xi)^{-1}\lesssim (\l^2\delta)^{-1}$ near
$Q_0^*\times \text{support}(\hat f)$.
Set
$$
\psi f=\sum_j\psi(x) q_j(D) f+\psi(x)\phi(D) f\,.
$$
Let
$q(x,\xi)=\psi(x)\phi(\xi)P_{\l\delta}(x,\xi)^{-1}$
smoothly localized in $\xi$ to $\{|\xi_1|\approx\l\,,|\xi_2|\le 2c\l\}\,.$
Then
$$
\psi(x)\phi(D) f=q(x,D)P_{\l\delta}f+r(x,D)f\,,
$$
where $r$ is of size $\l^{-1}\delta^{-2}$. Precisely, $q$ and $r$ are
supported where $|\xi|\approx\l$, and
$$
\l^2\delta\, q(\l^{-1}\delta^{-1}x,\l\,\delta\,\xi)\,,\quad
\l\,\delta^2 r(\l^{-1}\delta^{-1}x,\l\,\delta\,\xi)\in S^0_{0,0}\,.
$$
It follows that
$$
\|\psi(x)\phi(D) f\|_{H^1}\lesssim
\delta^{-\frac 32}|||f|||_{\l,\delta}\lesssim
\l^{\frac 14}|||f|||_{\l,\delta}\,.
$$
Since $p\ge 6$ then $\frac 14 \le \gamma(p+2)$, and Sobolev 
embedding yields
$$
\|\psi(x)\phi(D) f\|_{L^{p+2}}\lesssim
\l^{\gamma(p+2)}|||f|||_{\l,\delta}\,,
$$
hence we may take $\psi(x)\phi(D)f$ as the term $f_0$.

For each fixed $j$, consider the term $\psi(x)q_j(D)f$, and let $\xi_j$
be the center of the support of $q_j(\xi)$.
We can assume that the angular separation satisfies
$\angle(\xi_i,\xi_j)\gtrsim \delta\,|i-j|\,.$

Let $V_j$ denote the vector
$$
V_j = a_{\l\delta}(\xnot)\xi_j\,.
$$
Take a partition of unity in the $x_2$ variable, subordinate to a cover
by intervals of length $\delta^2$, and extend it to a partition of
unity in $(x_1,x_2)$ which is translation invariant under $V_j$, 
then multiply by $\psi$. We denote
the elements by $\psi_T(x):T\in \mathcal T_j$, so that
$\psi=\sum_{T\in\mathcal T_j}\psi_T$.
With $f_T=\psi_T q_j(D)f$, and $\mathcal T=\cup_j \mathcal T_j$, we have
the decomposition \eqref{decomp}.

Each $T=\text{support}(\psi_T)$ 
is contained in 
a tube of width $\delta^2$, length $\delta$, and direction $V_j$, 
where $V_j$ lies within a small angle of the $x_1$ axis.
We also have
\begin{equation}\label{psitbounds}
|(V_j\cdot d)^k\psi_T|\lesssim \l^k\delta^{-k}\,,\qquad\quad
|\partial_x^\alpha \psi_T|\lesssim \delta^{-2|\alpha|}\le 
\delta^{-2}(\l\delta)^{|\alpha|-1}\quad\text{if}
\quad |\alpha|\ge 1\,.
\end{equation}
We then expand $P_{\l\delta}(\psi_Tq_j(D)f)$ as
\begin{equation}\label{Ppsitdecomp}
(\dstar a_{\l\delta}d\psi_T)\, q_j(D)f
+2\langle a_{\l\delta}d\psi_T,d(q_j(D)f)\rangle
+\psi_T [P_{\l\delta},q_j(D)]f
+\psi_Tq_j(D)P_{\l\delta}f\,,
\end{equation}
and seek to show that
$$
\sum_j\sum_{T\in\mathcal T_j}
\l^{-2}\delta\,\|P_{\l\delta}(\psi_Tq_j(D)f)\|^2_{L^2}
\lesssim |||f|||^2_{\l,\delta}\,.
$$

For the first term in \eqref{Ppsitdecomp}, 
this follows by the finite overlap of the $\psi_T$
for $T\in \mathcal T_j$ and the finite overlap of the $q_j(\xi)$,
together with the pointwise bounds \eqref{psitbounds}.
The fourth term is similary handled by the finite overlap properties.

For the third term in \eqref{Ppsitdecomp}, 
we have the simple commutator bounds
$$
\|[P_{\l\delta},q_j(D)]\|_{L^2\rightarrow L^2}\le \l^2(\l\delta)^{-1}
=\l\delta^{-1}\,.
$$
Additionally, by the frequency localization of $a_{\l\delta}$ and 
$\rho_{\l\delta}$, the commutators have finite overlap as $j$ varies, 
yielding square summability over $j$.

We expand the brackets in the second term in \eqref{Ppsitdecomp} as
\begin{multline*}
\langle \,(a_{\l\delta}(x)-a_{\l\delta}(\xnot))
d\psi_T,dq_j(D)f\rangle
+i\langle V_j,d\psi_T\rangle\,q_j(D)f
\\
+\langle a_{\l\delta}(\xnot)
d\psi_T,(d-i\xi_j)(q_j(D)f)\rangle\,.
\end{multline*}
Each term has $L^2$ norm bounded by $\l\delta^{-1}\|f\|_{L^2}$, and
finite overlap properties yield square summability as above.

The finite overlap properties similarly yield that
$$
\sum_{T\in\mathcal T}\|f_T\|_{L^2}^2+\l^{-2}\|df_T\|_{L^2}^2\lesssim
\|f\|_{L^2}^2+\l^{-2}\|df\|_{L^2}^2\,,
$$
completing the verification of \eqref{sqsum}.\findemo

We also need a stronger inequality. For $T\in\mathcal T$, let $\xi_T$ equal
$\xi_j$ if $f_T=\psi_T(x) q_j(D)f$. Thus, the frequencies
of $f_T$ are concentrated in the $\l\delta$-ball about $\xi_T$, and
$|\xi_T|\approx\l$.

\begin{lemma}
The following bounds hold, for each $\alpha$, $\beta$,
\begin{equation}\label{sqsum2}
\left(\sum_T
\l^{-2|\alpha|}(\l\delta)^{-2|\beta|}
\|D^\alpha(D-\xi_T)^\beta f_T\|^2_{L^\infty_{x_1}L^2_{x_2}}
\right)^{\frac 12}
\le C_{\alpha,\beta}|||f|||_{\delta,\l}\,.
\end{equation}
\end{lemma}
\demo
Observe that we can write
$$
\l^{-|\alpha|}(\l\delta)^{-|\beta|}
D^\alpha(D-\xi_T)^\beta\psi_T(x)\,q_j(D)=\tilde\psi_T(x)\,\tilde q_j(D)
$$
where $\tilde\psi_T(x)$ and $\tilde q_j(\xi)$ satisfy similar support and
derivative bounds as $\psi_T$ and $q_j$. Hence, the proof we present for
the case $\alpha=\beta=0$ applies to the general case.

Let $q'_j(\xi)$ be a smooth
cutoff to the $\delta\l$-neighborhood of $\text{support} (q_j)$.
Then
$$
\|(1-q'_j(D))\psi_T q_j(D)f\|_{L^2}\lesssim\l^{-N}\|f\|_{L^2}\,.
$$
Since the number of tubes is bounded by $\delta^{-2}\ll\l$, Sobolev
embedding establishes the desired bounds on these terms.
We set $f'_T=q'_j(D)f_T$. By commutator arguments as above we
have $|||f'_T|||_{\delta,\l}\lesssim |||f_T|||_{\delta,\l}$.
The proof will then follow from \eqref{sqsum} by showing that
\begin{equation}\label{fluxbound}
\|f'_T\|_{L^\infty_{x_1}L^2_{x_2}}\lesssim |||f'_T|||_{\l,\delta}\,.
\end{equation}
We establish \eqref{fluxbound} by energy inequality arguments.
Let $V$ denote the vector field
$$
V=2(\partial_1\!f'_T)\,a_{\l\delta}\,df'_T+
\bigl(\l^2\rho_{\l\delta}\,{f'_T}^2-\langle a_{\l\delta}\,df'_T,df'_T
\rangle\bigr)
\overrightarrow{e_1}
$$
Then
$$
\dstar V=2(\partial_1\!f'_T)\,P_{\l\delta}f'_T
+\l^2(\partial_1\rho_{\l\delta}) {f'_T}^2-
\langle (\partial_1 a_{\l\delta})df'_T,df'_T\rangle
$$
Applying the divergence theorem on the set $x_1\le r$ yields
$$
\int_{x_1=r}V_1\,dx'
\lesssim
\l^2\|f'_T\|_{L^2}^2
+\delta^{-1}\|df'_T\|_{L^2}^2
+\delta\,\|P_{\l\delta}f'_T\|_{L^2}^2\le \l^2|||f'_T|||^2_{\l,\delta}
$$
Since $a_{\l\delta}$ and $\rho_{\l\delta}$ 
are pointwise close to the flat metric, we have pointwise that
$$
V_1\ge \tfrac 34 |\partial_1 f'_T|^2+\tfrac 34 \l^2|f'_T|^2
-\tfrac 32 |\partial_2f'_T|^2
$$
The frequency localization of $\widehat f'_T$ to $|\xi_2|\le c\l$ yields
$$
2\int_{x_1=r} V_1\,dx'\ge
\int_{x_1=r}|df'_T|^2+\l^2|f'_T|^2\,dx'\,.\qed
$$

\section{Overlap estimates}\label{section:overlap}
In this section we establish simple bounds on the overlap of tubes,
and resulting $\ell^p$ bounds on the counting function.

\begin{lemma}\label{overlap}
Let $x$ and $y$ be two points in $Q_0$. Then the number of
distinct tubes $T\in \mathcal T$ which pass within distance $4\delta^2$
of both $x$ and $y$ is bounded by
$C\min\bigl(\delta^{-1},\frac{\delta}{|x_1-y_1|}\bigr)$.
\end{lemma}
\demo
For each $j$, there is a fixed bound on the number of
tubes $T\in\mathcal T_j$ which pass within distance $4\delta^2$ of $x$.
It thus suffices to bound the number of distinct $j$ such that
the line through $x$ in direction $V_j$
passes within distant $\sim\delta^2$ of $y$.
The above bound is then a simple consequence
of the fact that $\angle(V_i,V_j)\gtrsim \delta|\,i-j|$.
\qed

\medskip

Now consider a collection $\mathcal N\subset \mathcal T$ 
containing $N$ distinct tubes. We make a
decomposition of the cube $Q_0$ into a $\delta^{-1}\times\delta^{-1}$
grid $\mathcal Q$
of cubes $Q$ of sidelength $\delta^2$. Let $n_Q$ denote the number
of tubes in $\mathcal N$ which intersect $Q^*$,
$$
n_Q=\#\{T\in\mathcal N:T\cap Q^*\ne\emptyset\}\,.
$$
By $\|n_Q\|_{\ell^p\ell^q}$ we mean the $\ell^p_{x_1}\ell^q_{x_2}(\mathcal Q)$
norm of the counting function $n_Q$, taken over the 
grid $\mathcal Q$ of cubes.

\begin{corollary}
The following bounds hold,
$$
\|n_Q\|_{\ell^\infty\ell^1}\lesssim N\,,
\qquad\qquad
\|n_Q\|_{\ell^2\ell^\infty}\lesssim 
|\log\delta|^{\frac 12}\delta^{-\frac 12}N^{\frac 12}\,.
$$
Furthermore, for $q\ge 3$
\begin{equation}\label{tubebound}
\|n_Q\|_{\ell^q(\mathcal Q)}\lesssim
|\log\delta|^{\frac 1q}\delta^{-\frac 1q}N^{1-\frac 1q}\,.
\end{equation}
\end{corollary}
\demo
The first bound is an immediate consequence of the fact that, for each
$T$ and each $r$, 
there is a fixed upper bound on the number of cubes $Q^*$ centered
on the line $x_1=r$ that intersect $T$. For the second bound, 
we consider the map
$$
%\{c_T\}\in \ell^2(\mathcal T)\rightarrow
W\{c_T\}=\sum_T c_T\chi_T(Q)\quad\text{where}\quad
\begin{cases}\chi_T(Q)=1\,,& T\cap Q^*\ne \emptyset\\
\chi_T(Q)=0\,,& T\cap Q^*= \emptyset\end{cases}
$$
It suffices to show that 
$W:\ell^2(\mathcal T)\rightarrow \ell^2\ell^\infty(\mathcal Q)$
with bound $|\log\delta|^\frac 12\delta^{-\frac 12}$. The map $WW^*$ takes
the form
$$
WW^*\{c_Q\}(Q')=\sum_{Q}n(Q',Q)c_Q\,,\qquad
n(Q',Q)=\#\{\,T:T\cap Q^*\ne\emptyset\,, T\cap {Q'}^*\ne\emptyset\,\}\,.
$$
We need to show that 
$WW^*:\ell^2\ell^1(\mathcal Q)\rightarrow \ell^2\ell^\infty(\mathcal Q)$
with norm $|\log\delta|\,\delta^{-1}$. This is an easy consequence
of the bound from Lemma \ref{overlap},
$$
n(Q',Q)\lesssim \min\Bigl(\delta^{-1},\,\delta\,|x_1(Q')-x_1(Q)|^{-1}\Bigr)\,.
$$
Applying interpolation now yields the bounds
$$
\|n_Q\|_{\ell^q\ell^r}\lesssim 
|\log\delta|^{\frac 1q}\delta^{-\frac 1q}N^{1-\frac 1q}
\,,
\qquad \tfrac 2q +\tfrac 1r = 1\,. 
$$
Note that if $q\ge 3$ then $r\le q$, yielding \eqref{tubebound}.
\findemo

\section{Proof of Theorem \ref{maintheorem'}}\label{section:proof}
Given $f$ and the cube $Q_0$, we decompose 
$\psi f=\sum_{T\in \mathcal T}f_T+f_0$ as in
Section \ref{section:tube}, and control $\|f_0\|_{L^{p+2}}$ 
by Sobolev embedding.
We make a further decomposition by collecting together tubes for which
$f_T$ is of comparable energy. Precisely, decompose 
$\mathcal T=\cup_{k\ge -k_0}\,\mathcal N_k$ where
%, with $C$ the constant from \eqref{sqsum}, then 
$T\in\mathcal N_k$ if
\begin{multline*}
2^{-k-1}|||f|||_{\l,\delta}\;<\;|||f_T|||_{\l,\delta}\,+\!\!
\sum_{|\alpha|+|\beta|\le 3}\l^{-|\alpha|}(\l\delta)^{-|\beta|}
\|D^{\alpha}(D-\xi_T)^\beta f_T\|_{L^\infty_{x_1}L^2_{x_2}}\\
\le 2^{-k}|||f|||_{\l,\delta}\,.
\end{multline*}
We handle the tubes for $k\ge 2\log_2\!\l$ by the Sobolev bound 
$\|f_T\|_{L^\infty}\le \|Df_T\|_{L^\infty_{x_1}L^2_{x_2}}$,
since there are at most $\l^{\frac 13}$ tubes in all.
This leaves at most $\approx\log\l$ values of $k$, 
which we handle individually.
We thus fix some $\mathcal N=\mathcal N_k$, and let $N$ be the number
of tubes in $\mathcal N$. We multiply $f$ by a constant so that
$|||f|||_{\l,\delta}= 2^k$, which by \eqref{sqsum} and \eqref{sqsum2}
implies $N^{\frac 12}\lesssim |||f|||_{\l,\delta}$.
We then need to establish the following.

\begin{theorem}\label{maintheorem''}
Suppose that $f=\sum_{T\in\mathcal N}f_T$, where each $f_T$ is supported
in $Q_0$, and 
$$
|||f_T|||_{\delta,\l}+
\sum_{|\alpha|+|\beta|\le 3}\l^{-|\alpha|}(\l\delta)^{-|\beta|}
\|D^{\alpha}(D-\xi_T)^\beta f_T\|_{L^\infty_{x_1}L^2_{x_2}}
\le 1\,.
$$
Let $N$ denote the cardinality of $\mathcal N$. Then, under the conditions
of Hypothesis \ref{hypo},
$$
\|f\|_{L^{p+2}}\lesssim (\log\l)^{p-5}\l^{\gamma(p+2)}N^{\frac 12}\,.
$$
\end{theorem}
\demo
As above we decompose $Q_0$ into cubes $Q$ of size $\delta^2$,
$Q_0=\cup_{\mathcal Q}Q$. 
By hypothesis, for each $Q$ we have
$$
\|f\|_{L^p(Q)}\lesssim (\log\l)^{p-6}\l^{\gamma(p)}
|||\sum_{T\cap Q^*\ne\emptyset}f_T\,|||_{\delta^2,\l,Q^*}\,.
$$
We first show that
\begin{equation}\label{tsum}
|||\sum_{T\cap Q^*\ne\emptyset}f_T\,|||_{\delta^2,\l,Q^*}
\lesssim n_Q^{\frac 12}\,.
\end{equation}
For this, note that $|a_{\l\delta}-a_{\l\delta^2}|\le(\l\delta^2)^{-1}$,
hence
$$
\bigl|\bigl|\bigl|\sum_{T\cap Q^*\ne\emptyset}f_T\,
\bigr|\bigr|\bigr|_{\delta^2,\l,Q^*}\le
\sum_{|\alpha|\le 2}\l^{-|\alpha|}
\bigl\|\sum_{T\cap Q^*\ne\emptyset}D^\alpha\! f_T\,
\bigr\|_{L^\infty_{x_1}L^2_{x_2}}
+\l^{-1}\delta\,
\bigl\|\sum_{T\cap Q^*\ne\emptyset}P_{\l\delta}f_T\,\bigr\|_{L^2}\,.
$$
For each $j$, there are a bounded number of tubes $T$ with $\xi_T=\xi_j$ 
for which $T\cap Q^*\ne\emptyset\,,$ hence we can assume
the different $\xi_T$ in the above sum are spaced by distance $\l\delta$
in the $\xi_2$ variable. Thus,
\begin{align*}
\l^{-|\alpha|}
\bigl\|\sum_{T\cap Q^*\ne\emptyset}D^\alpha\!f_T
\bigr\|_{L^\infty_{x_1}L^2_{x_2}}
&\lesssim
\left(\sum_{|\beta|\le 1}\,
\sum_{T\cap Q^*\ne\emptyset}
\l^{-2|\alpha|}(\l\delta)^{-2|\beta|}
\|(D-\xi_T)^\beta D^\alpha\! f_T\|^2_{L^\infty_{x_1}L^2_{x_2}}
\right)^{\frac 12}\\
&\lesssim n_Q^{\frac 12}\,.
\end{align*}
To complete the proof of \eqref{tsum}, we use that $n_Q\le\delta^{-1}$
to bound
$$
\l^{-1}\delta\sum_{T\cap Q^*\ne\emptyset}\|P_{\l\delta}f_T\|_{L^2}
\le
\delta^{\frac 12}\!\!
\sum_{T\cap Q^*\ne\emptyset}|||f_T|||_{\delta,\l}
\le \delta^{\frac 12}n_Q
\le n_Q^{\frac 12}\,.
$$
By  \eqref{tsum} and Hypothesis \ref{hypo}, we thus have
\begin{equation}\label{ftpbound}
||f||_{L^p(Q)}
\lesssim (\log\l)^{p-6}\l^{\gamma(p)}n_Q^{\frac 12}\,.
\end{equation}
We next note the bound
\begin{align*}
\|f_T\|_{L^\infty}&\le
2\|(D-\xi_T)f_T\|_{L^\infty_{x_1}L^2_{x_2}}^{\frac 12}
\|f_T\|_{L^\infty_{x_1}L^2_{x_2}}^{\frac 12}\\
&\le \l^{\frac 12}\delta^{\frac 12}
\sum_{|\beta|\le 1}(\l\delta)^{-|\beta|}
\|(D-\xi_T)^\beta f_T\|_{L^\infty_{x_1}L^2_{x_2}}\\
&\le \l^{\frac 12}\delta^{\frac 12}\,.
\end{align*}
Consequently,
\begin{equation}\label{ftinftybound}
\|f\|_{L^\infty(Q)}\le \l^{\frac 12}\delta^{\frac 12}n_Q\,.
\end{equation}
%%%%%New version, old one below.
Combining \eqref{ftpbound}--\eqref{ftinftybound} with \eqref{tubebound}
for $q=\frac 12 p$ and $q=\infty$ respectively, we obtain
\begin{align*}
\|f\|_{L^p} & \lesssim 
(\log\l)^{p-6+\frac 1p}\l^{\gamma(p)}\delta^{-\frac 1p}N^{\frac 12-\frac 1p}\\
\|f\|_{L^\infty} & \lesssim \l^{\frac 12}\delta^{\frac 12}N
\end{align*}
Interpolation yields
$$
\|f\|_{L^{p+2}}\lesssim 
(\log\l)^{\frac{p(p-6)+1}{p+2}}\l^{\gamma(p+2)}N^{\frac 12}
$$
Observing that if $p\ge 6$ we have $p(p-6)+1\le (p+2)(p-5)$ concludes
the proof.
\qed


\begin{thebibliography}{10}

\bibliographystyle{plain}

\bibitem{GeTa} D. Geba and D. Tataru,
\newblock{Dispersive estimates for wave equations},
Comm. Partial Differential Equations {\bf 30} (2005), no 4-6, 849--880.

\bibitem{Gr} D. Grieser,
\newblock{\em $L^p$ bounds for eigenfunctions and spectral
projections of the Laplacian near concave boundaries}, Thesis, UCLA, 1992.

\bibitem{KST} H. Koch, H. Smith and D. Tataru,
\newblock{\em Sharp $L^q$ bounds on spectral clusters
for Holder metrics},
\newblock Math. Res. Lett. {\bf 14} (2007), no. 1, 77--85.

\bibitem{KTZ} H. Koch, D. Tataru and M. Zworski,
\newblock{\em Semiclassical $L^p$ estimates}, 
to appear in Annales Henri Poincar\'e.

\bibitem{Sm1} H. Smith,
\newblock{\em Spectral cluster estimates for $C^{1,1}$ metrics},
\newblock Amer. Jour. Math. {\bf 128} (2006), 1069--1103.

\bibitem{Sm2} H. Smith,
\newblock{\em Sharp $L^2 \rightarrow L^q$ bounds on spectral projectors
for low regularity metrics},
\newblock Math. Res. Lett. {\bf 13} (2006), no. 6, 967--974.

\bibitem{SmSo1} H. Smith and C. Sogge,
\newblock{\em On Strichartz and eigenfunction estimates for low regularity
metrics},
\newblock Math. Res. Lett. {\bf 1} (1994), 729--737.

\bibitem{SmSo2} H. Smith and C. Sogge,
\newblock{\em On the $L^p$ norm of spectral clusters for compact 
manifolds with boundary},
\newblock Acta Math. {\bf 198} (2007), 107--153.

\bibitem{So} C. Sogge,
\newblock{\em Concerning the $L^p$ norm of spectral clusters for second
order elliptic operators on compact manifolds},
\newblock J. Funct. Anal. {\bf 77} (1988), 123--134.

\end{thebibliography}
\end{document}